\documentclass[final,twocolumn]{elsarticle}

\usepackage{amssymb}
\usepackage{amsmath}

\begin{document}

\begin{frontmatter}



\title{Two Methods for Numerical Inversion of the Z-Transform}


\author{Farshad Merrikh-Bayat}

\address{Department of Electrical and Computer Engineering, University of Zanjan, Zanjan, Iran. Email:
f.bayat@znu.ac.ir}

\begin{abstract}
In some of the problems, complicated functions of the Z-transform
variable, $z$, appear which either cannot be inverted analytically
or the required calculations are quite tedious. In such cases
numerical methods should be used to find the inverse Z-transform.
The aim of this paper is to propose two simple and effective
methods for this purpose. The only restriction on the signal
(whose Z-transform is given) is that it must be absolutely
summable (of course, this limitation can be removed by a suitable
scaling). The first proposed method is based on the Discrete
Fourier Transform (DFT) and the second one is based on solving a
linear system of algebraic equations, which is obtained after
truncating the signal whose Z-transform is known. Numerical
examples are also presented to confirm the efficiency of the
proposed methods. Functions in non-integer powers of $z$ are also
briefly discussed and it is shown that such functions cannot be
obtained by taking the Z-transform from any discrete-time signal.
\end{abstract}


\end{frontmatter}


\section{Introduction}
In some of the problems we may face with complicated functions of
the Z-transform variable, $z$, whose inverse transform either
cannot be calculated analytically or the required calculations are
quite tedious\footnote{Please note that this is a draft of a paper
whose extended version may be submitted to a conference or journal
for possible publication. The last version may also be subjected
to modifications during the review process.}. Clearly, numerical
methods should be used in dealing with such problems to find the
inverse Z-transform. At this time, as far as the author knows, the
main contribution in the field of numerical inversion of
Z-transform is \cite{pop}, which may seem complicated for software
implementation. Moreover, so far no application of this method is
reported in the literature. The aim of this brief paper is to
represent two very simple and effective methods for numerical
inversion of the Z-transform, which can be easily implemented by
software.

As a well-known classical fact, $x[n]$ and its Z-transform $X(z)$
are related through the following equation
\begin{equation}\label{int_z}
x[n]=\frac{1}{2\pi j}\oint_c X(z) z^{n-1}dz,
\end{equation}
where $c$ is a closed simple counterclockwise contour considered
in the region of convergence (ROC). At the first looking it may
seem that the most straightforward method for obtaining the
inverse Z-transform of the given $X(z)$ is to calculate the
integral in (\ref{int_z}) numerically. For this purpose, assuming
that the ROC of $X(z)$ is equal to $|z|>a$ $(a>0)$ we can perform
the numerical integration on the contour $c:re^{j\theta}$, where
$0<\theta<2\pi$ and $r$ is an arbitrarily chosen constant such
that $r>a$. For the absolutely summable $x[n]$, $r$ can be
considered equal to unity which leads to
\begin{equation}\label{int_four}
x[n]=\frac{1}{2\pi}\int_0^{2\pi} X(e^{j\omega}) e^{j\omega n}
d\omega.
\end{equation}
The above equation simply represents the inverse Fourier transform
of $X(e^{j\omega})=X(z)|_{z=\exp(j \omega)}$.

The main problem with this approach (as well as all of the others
based on the numerical evaluation of an integral) is that
sometimes it does not lead to accurate results. More precisely,
according to the Residue theorem the $x[n]$ calculated from
(\ref{int_z}) must be independent of the special contour used for
integration (which, of course, satisfies the general conditions).
But in practice it is observed that when the contour is considered
as a circle inside the ROC, changing the radius of circle
considerably changes the amplitude of the samples of $x[n]$ (in
the numerical examples presented at the end of this paper it is
observed that the most accurate results are obtained when the
radius of circle is considered equal to unity). Especially, the
numerical integration algorithms may lead to less accurate results
for larger values of $n$ (note that the frequency of the kernel of
(\ref{int_four}) is increased by increasing $n$). The other
problem with this technique is that for any value of $n$ we have
to evaluate the integral again, which may be time consuming. By
the way, calculation of the inverse Z-transform by numerical
evaluation of an integral is a well-known method widely used by
researchers. The aim of this paper is to discuss on two less
famous methods.

In the following discussions it is assumed that the signal $x[n]$
whose Z-transform is known is equal to zero for all $n<0$.
Moreover, it is assumed that the signal under consideration is
absolutely summable which is equal to the fact that all poles of
$X(z)$ are located inside the unit circle.\\

\textbf{Remark 1.} The methods proposed in this paper can still be
used even if $x[n]$ is not absolutely summable. For this purpose
we can simply use the equation $F\left\{a^nx[n]\right\}
=X(a^{-1}z)$ for some $|a|<1$ to map all poles of $X(z)$ to the
region inside the unit circle.\\

The rest of this paper is organized as the following. The first
and the second method for numerical inversion of the Z-transform
are explained in Sections 2 and 3, respectively. Two numerical
examples are studied in Section 4.

\section{The first method for numerical inversion of the Z-transform}
The first method is based on truncating the absolutely summable
signal whose Z-transform is known, forming a linear system of
algebraic equations, and then solving it. According to the
definition of Z-transform, $x[n]$ and $X(z)=Z\{ x[n]\}$ are
related through the following equation:
\begin{equation}\label{sec1}
X(z)=\sum_{n=0}^\infty x[n] z^{-n}.
\end{equation}
(Recall that, without a considerable loss of generality, in this
paper it is assumed that $X(z)$ does not have any poles at
infinity, i.e. $\lim_{|z|\rightarrow\infty}X(z) <\infty$, which is
equivalent to the fact that $x[n]=0$ for $n<0$). Clearly, Eq.
(\ref{sec1}) holds for any $z$ in the ROC. For example, for
$z_1\in\mathrm{ROC}$ we have
\begin{equation}\label{sec2}
X(z_1)=\sum_{n=0}^\infty x[n] z_1^{-n}.
\end{equation}
Considering the fact that $x[n]$ is absolutely summable, in the
above equation $x[n]$ tends to zero as $n$ tends to infinity.
Moreover, the weights $z_1^{-n}$ in (\ref{sec2}) also rapidly tend
to zero by increasing the value of $n$ provided that $|z_1|>1$
(note that according to the above discussion the region outside
the unit circle necessarily belongs to ROC). It concludes that the
sigma in the righthand side of (\ref{sec2}) can safely be
approximated by its first few terms as the following:
\begin{equation}\label{sec3}
X(z_1)\approx\sum_{n=0}^N x[n] z_1^{-n},
\end{equation}
where $N$ is chosen sufficiently large.

Assuming that $z_1, z_2, \ldots, z_m$ are (random) points chosen
from ROC, Eq. (\ref{sec3}) leads to the following approximate
linear algebraic equation:
\begin{multline}\label{eq1}
\left(%
\begin{array}{c}
  X(z_1) \\
  X(z_2) \\
  \vdots \\
  X(z_m) \\
\end{array}%
\right)\approx\left(%
\begin{array}{ccccc}
  1 & z_1^{-1} & z_1^{-2} & \ldots & z_1^{-N} \\
  1 & z_2^{-1} & z_2^{-2} & \ldots & z_2^{-N} \\
  \vdots & ~ & \ddots & ~ & \vdots \\
  1 & z_m^{-1} & z_m^{-2} & \ldots & z_m^{-N} \\
\end{array}%
\right)\\ \times \left(%
\begin{array}{c}
  x[0] \\
  x[1] \\
  \vdots \\
  x[N] \\
\end{array}%
\right),
\end{multline}
where $m$ is the number of (random) points chosen from ROC.
Assuming $\mathbf{X}=[X(z_1), \ldots, X(z_m)]^T$,
$\mathbf{x}=[x[0], \ldots, x[N]]^T$, and
\begin{equation}\label{a}
\mathbf{A}=\left(%
\begin{array}{ccccc}
  1 & z_1^{-1} & z_1^{-2} & \ldots & z_1^{-N} \\
  1 & z_2^{-1} & z_2^{-2} & \ldots & z_2^{-N} \\
  \vdots & ~ & \ddots & ~ & \vdots \\
  1 & z_m^{-1} & z_m^{-2} & \ldots & z_m^{-N} \\
\end{array}%
\right),
\end{equation}
Eq. (\ref{eq1}) can be written as $\mathbf{Ax}=\mathbf{X}$ which
contains $N$ variables and $m$ equations. Assuming $m=N$ this
equation (theoretically) has a unique solution provided that the
points $z_i$ ($i=1,\ldots,m$) are suitably chosen such that the
$\mathbf{A}$ given in (\ref{a}) is full rank.

In the simulations of this paper the value of $m$ is considered
slightly larger than $N$ (more precisely, $m\approx 1.1N$) and the
points $z_i$ ($i=1,\ldots,m$) are randomly chosen from the region
defied by $1<|z|<2$. Then the solution of (\ref{eq1}) is obtained
by finding the $\mathbf{x}$ which minimizes $||\mathbf{Ax-X}||_2$
(such a solution can easily be obtained by using the Matlab
command '$\backslash$', i.e., $\mathbf{x=A\backslash X}$). \\

\textbf{Remark 2.} In practice it is observed that using the $z_i$
such that $|z_i|<1$ makes the problem ill-conditioned for larger
values of $N$ (even if $z_i$ belongs to ROC). Hence, points inside
the unit circle should be avoided.

\section{The second method for numerical inversion of the Z-transform}
\subsection{Mathematical background}
As mentioned earlier, the second proposed method is based on the
discrete Fourier transform (DFT), which can effectively be
calculated through the fast Fourier transform (FFT) (see the
\textsf{fft} and \textsf{ifft} commands in Matlab). Before
explaining our method, we need to briefly review the main result
of DFT. The following discussion can be found in \cite{oppenheim2}
with more details.

Consider the discrete-time signal $x[n]$ for which we have
$x[n]=0$, $n<0$, and $F\{x[n]\}=X(e^{j\omega})$. Obviously
$X(e^{j\omega})$ is periodic with period $2\pi$. Let us take $N$
samples of $X(e^{j\omega})$ at each period and denote it as
$\widetilde{X}[k]$; that is $\widetilde{X}[k]=
X(e^{j\omega})|_{\omega=2\pi k/N}$, $k\in \mathbb{Z}$. It can be
easily verified that $\widetilde{X}[k]$ is also periodic with
period $N$. Now, one can prove that \cite{oppenheim2} if
$\widetilde{X}[k]$ ($k=0,\ldots,N-1$) is considered as the
coefficients of a periodic signal like $\widetilde{x}[n]$, then
$x[n]$ and $\widetilde{x}[n]$ are related through the following
equation:
\begin{equation}\label{dft1}
\widetilde{x}[n]=\sum_{r=-\infty} ^\infty x[n-rN].
\end{equation}
Considering the fact that $x[n]=0$ for $n<0$, Eq. (\ref{dft1}) can
further be written as
\begin{equation}\label{dft2}
\widetilde{x}[n]=\sum_{r=-\infty} ^0 x[n-rN].
\end{equation}

According to the above discussion, the inverse discrete Fourier
transform (IDFT) of $\widetilde{X}[k]$ (i.e., the IDFT of samples
of $ X(e^{j\omega})=F\{x[n]\}$) is equal to $\widetilde{x}[n]$ as
given in (\ref{dft2}). But, $\widetilde{x}[n]$ is an approximation
for $x[n]$ provided that $\lim_{n\rightarrow\infty}x[n]=0$ (which
is the case in our problem since $x[n]$ is assumed to be
absolutely summable) and $N$ is chosen sufficiently large.

It concludes that in order to find the inverse Z-transform of
$X(z)$, whose all poles are located inside the unit circle and has
no poles at infinity, first we calculate $X(e^{j\omega})
=X(z)|_{z=\exp(j\omega)}$ (which is valid since the unit circle
belongs to ROC according to our previous assumption). Then we
calculate the samples of $X(e^{j\omega})$ at frequencies
$\omega=2\pi k/N$, where $k=0,\ldots,N-1$ and $N$ is a
sufficiently big number. Denote the $k$th sample of
$X(e^{j\omega})$ as $X[k]$. Taking the IDFT of $X[k]$ (using,
e.g., the Matlab command \textsf{ifft}) we arrive at
$\widetilde{x}[n]$ whose first $N$ samples (i.e.,
$\widetilde{x}[n]$ for $n=0,\ldots,N-1$) are approximately equal
to $x[n]$. The reason for this statement is that according to
(\ref{dft2}) we have
\begin{equation}\label{dft3}
\widetilde{x}[n]=x[n]+\sum_{r=-\infty} ^{-1} x[n-rN].
\end{equation}
But the sigma in the righthand side of (\ref{dft3}) converges to
zeros very rapidly and consequently we have $\widetilde{x}[n]
\approx x[n]$ (more detailed discussion can be found in the
following).

\subsection{Comprehensive algorithm for numerical inversion of the Z-transform (second method)}
The following is a brief explanation of the second method for
calculation of the inverse Z-transform of $X(z)$ whose all poles
are located inside the unit circle and has no poles at infinity:

\begin{enumerate}
    \item Calculate $X(e^{j\omega})=X(z)|_{z=\exp(j\omega)}$.
    \item Calculate $X[k]=X(e^{j\omega})|_{\omega=2\pi k/N}$,
    $k=0, \ldots, N-1$ (Increasing the value of $N$ increases the accuracy of results).
    \item Calculate the IDFT of $X[k]$ (using the Matlab command
    \textsf{ifft}) to arrive an approximation like $\widetilde{x}[n]$ for $x[n]=Z^{-1}\{ X(z)\}$.
\end{enumerate}

\subsection{Error analysis}
In the following we discuss on the error caused by the above
method in a simple case. According to (\ref{dft3}) we have
\begin{equation}\label{dft4}
\widetilde{x}[n]=x[n]+\sum^\infty_{r=1} x[n+rN].
\end{equation}
The second term in the righthand side of (\ref{dft4}) indicates
the error caused by the proposed method, i.e., if
$\sum^\infty_{r=1} x[n+rN]=0$ then $\widetilde{x}[n]=x[n]$ and the
inversion is errorless.

Assuming that $X(z)=F\{x[n]\}$ is a rational function in $z$ whose
farthest pole from the origin is simple and located at $z=a$
($|a|<1$), for sufficiently large values of $n$ we have
\begin{equation}\label{dft5}
x[n]\sim A a^n\quad (n\rightarrow\infty),
\end{equation}
for some $A\in\mathbb{R}$. Hence, the corresponding error in the
$n$th sample of $x[n]$ (assuming sufficiently large values for
$N$) is obtained as the following:
\begin{equation}
\mathrm{error~ in}~ \widetilde{x}[n] =\sum^\infty_{r=1}
x[n+rN]\approx \sum^\infty_{r=1} A a^{n+rN},
\end{equation}
which leads to:
\begin{equation}\label{err_dft}
\mathrm{absolute ~error~ in}~
\widetilde{x}[n]\approx\frac{Aa^{n+N}}{1-a^N}, \quad n=0,\ldots,
N-1.
\end{equation}
(Note that according to (\ref{dft4}) the error in, e.g.,
$\widetilde{x}[0]$ is equal to $x[N]+x[2N]+x[3N]+\ldots$. On the
other hand, for sufficiently large values of $N$ the approximation
given in (\ref{dft5}) is valid which leads to the error calculated
above for the $n$th sample of $x[n]$).

It can be easily verified that the error calculated from
(\ref{err_dft}) is a monotonically increasing function of $a$,
that is, smaller the value of $|a|$ (which is equal to the
amplitude of the largest pole of $X(z)$) smaller the amplitude of
error in all samples of $\widetilde{x}[n]$. Moreover, it is
observed that increasing the value of $N$ is highly effective for
decreasing the error in $\widetilde{x}[n]$. In fact, for a fixed
$n$ the error in $\widetilde{x}[n]$ is approximately proportional
to $a^N$, which means that duplication of $N$ decreases the error
in $\widetilde{x}[n]$ by the factor of $a^N$ (for example,
assuming $|a|=0.9$, application of $N=40$ instead of $N=20$
decreases the error $0.9^{-20}\approx 8.2$ times). Fortunately,
the FFT algorithm is capable of working with large values of $N$.
Finally, note that according to (\ref{err_dft}) the absolute error
is larger in the first samples of $\widetilde{x}[n]$ (i.e.,
smaller values of $n$).

A similar discussion can be presented for rational $X(z)$ with
repeated poles, which is not discussed here.\\

\textbf{Remark 3.} Note that a function in non-integer
(fractional) powers of $z$ cannot be obtained by taking the
Z-transform from any real-world discrete-time signal. For example,
there is no $x[n]$ such that $X(z)=Z\{x[n]\}=1/(1-0.5z^{-0.5})$
(however, there is a $x(t)$ such that, e.g., $X(s)=L\{x(t)\}=1
/(s^{0.5}+1)$, where $L$ stands for the Laplace transform). The
reason is that $X(z)$ is necessarily single-valued by definition
which cannot be equal to any multi-valued (fractional-order)
function of $z$. That is why the fractional calculus does not have
a discrete-time dual.

\section{Numerical Examples}
Two numerical examples are presented in this section in support of
the proposed methods. In each case, the numerical inversion is
performed directly by using (\ref{int_four}) as well as the
proposed algorithms, and then the results are compared. In all of
the following simulations, the integrals in (\ref{int_z}) and
(\ref{int_four}) are evaluated numerically using the Matlab R2009a
command \textsf{quadgk}.\\

\begin{figure}[tb]
\begin{center}
\includegraphics[width=8.5cm]{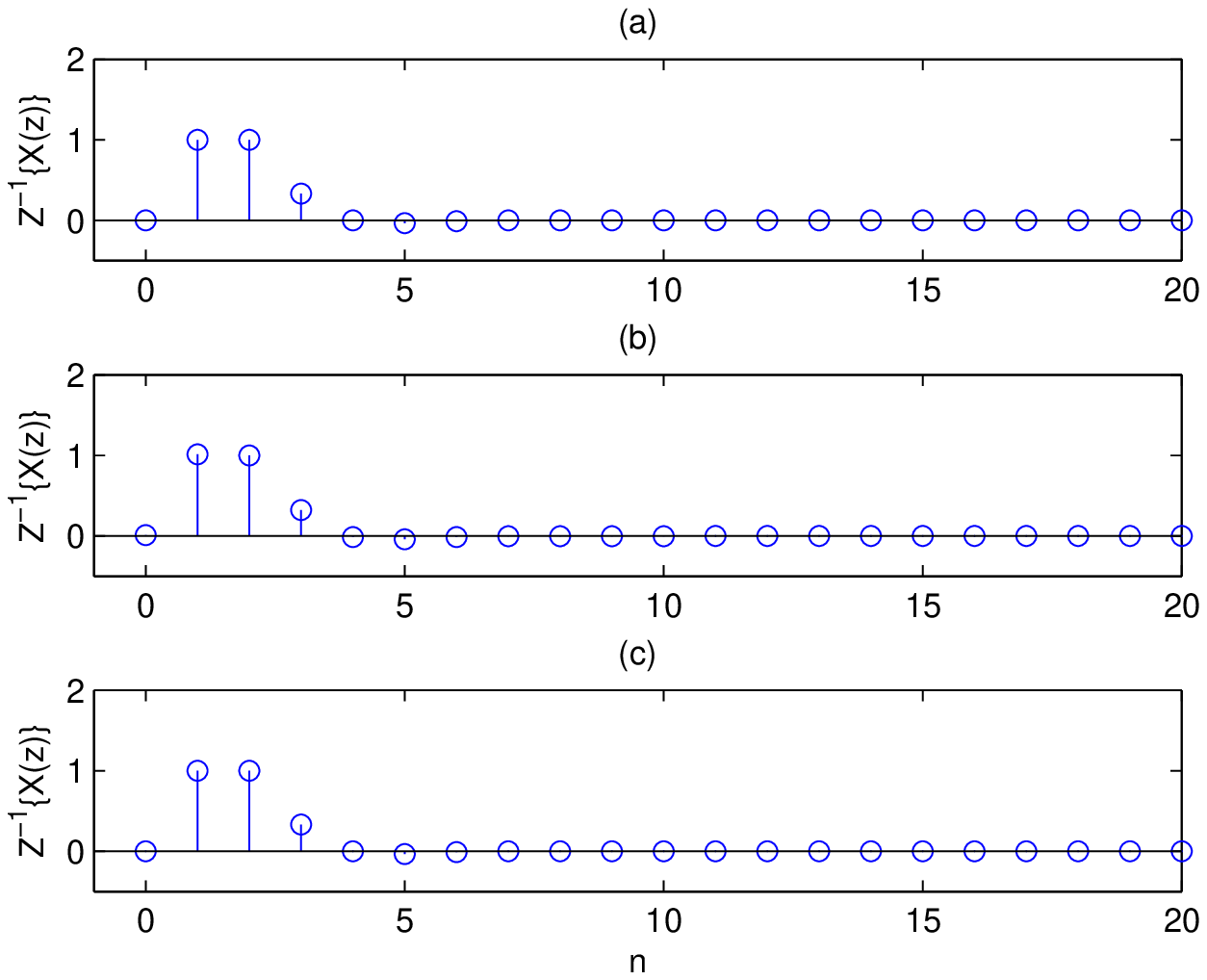}
\caption{Inverse Z-transform of $X(z)$ when (a) the first method,
(b) the second method, and (c) Eq. (\ref{int_four}) is applied
(corresponding to Example 1).}\label{fig_ex11}
\end{center}
\end{figure}

\begin{figure}[tb]
\begin{center}
\includegraphics[width=8.5cm]{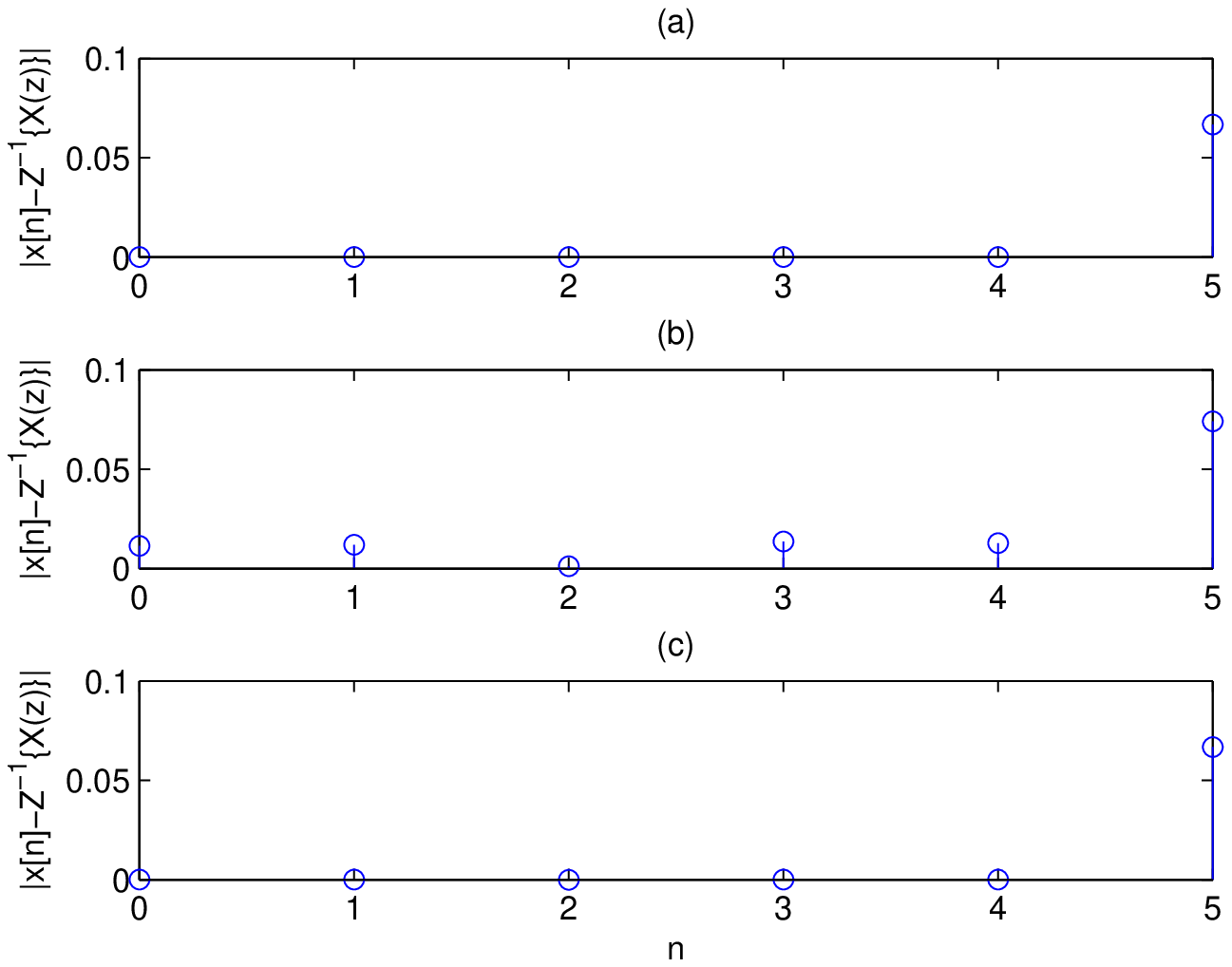}
\caption{The absolute error caused by (a) the first method, (b)
the second method, and (c) Eq. (\ref{int_four}), in the fist six
samples of $x[n]$ (corresponding to Example 1).}\label{fig_ex12}
\end{center}
\end{figure}

\textbf{Example 1.} Consider $X(z)=e^\frac{1}{z} \sin
\left(\frac{1}{z}\right)$ whose ROC is equal to $|z|>0$. The
inverse Z-transform of this function can be obtained through the
Laurent series expansion as the following
\begin{align*}
X(z)&=e^\frac{1}{z} \sin \left(\frac{1}{z}\right)\\
&=\left(1+\frac{1}{z}+\frac{1}{2!z^2}+\ldots\right)\\
&\times\left(\frac{1}{z}-\frac{1}{3!z^3}+\frac{1}{5!z^5}-\ldots\right)\\
&=\frac{1}{z}+\frac{1}{z^2}+\left(-\frac{1}{3!}+\frac{1}{2!}\right)
\frac{1}{z^3}+ \left(-\frac{1}{3!}+\frac{1}{3!}\right)
\frac{1}{z^4}\\
&+\left(\frac{1}{5!}-\frac{1}{2!\times 3!}+\frac{1}{4!}\right)
\frac{1}{z^5}+\ldots\\
&=z^{-1}+z^{-2}+\frac{1}{3}z^{-3}-\frac{1}{30}z^{-5}+\ldots
\end{align*}
which concludes that
\begin{equation}
x[0]=0,\quad x[1]=x[2]=1, \quad x[3]=\frac{1}{3}
\end{equation}
\begin{equation}
x[4]=0,\quad x[5]=\frac{1}{30}, \ldots
\end{equation}
and $x[n]=0$, $n<0$. Figs. \ref{fig_ex11}(a)-(c) show the inverse
Z-transform of $X(z)$ when the first method, second method, and
Eq. (\ref{int_four}) is applied, respectively. Figs.
\ref{fig_ex12}(a)-(c) show the absolute error caused by each of
these methods, respectively.

\begin{figure}[tb]
\begin{center}
\includegraphics[width=8.5cm]{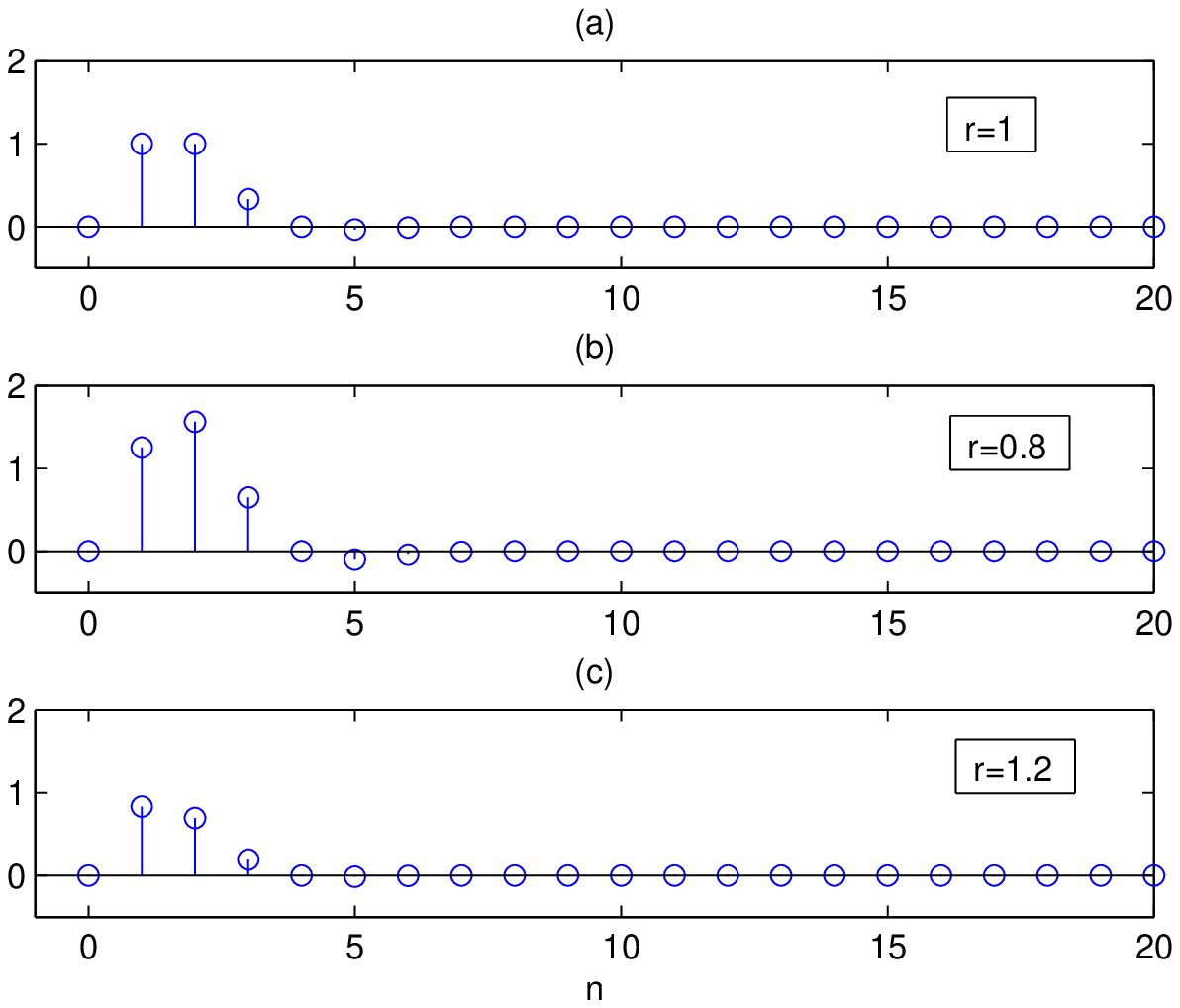}
\caption{The $x[n]$ obtained through (\ref{int_z}) when $c$ is
considered as a circle with radius $r=1,0.8,1.2$ (corresponding to
Example 1). $r=1$ leads to the most accurate
results.}\label{fig_ex13}
\end{center}
\end{figure}

As mentioned earlier, in case of using (\ref{int_z}) we will see
that (in practice) the value of $x[n]$ depends on the especial
shape of the contour $c$. Figs. \ref{fig_ex13}(a)-(c) show the
$x[n]$ obtained through (\ref{int_z}) when $c$ is considered as a
circle with radius $r=1,0.8,1.2$, respectively. As it can be
observed in this figure, the radius of contour highly affects the
resulted $x[n]$. However, comparing Fig. \ref{fig_ex13} with Fig.
\ref{fig_ex11} concludes that the most accurate result corresponds
to $r=1$.\\

\begin{figure}[tb]
\begin{center}
\includegraphics[width=8.5cm]{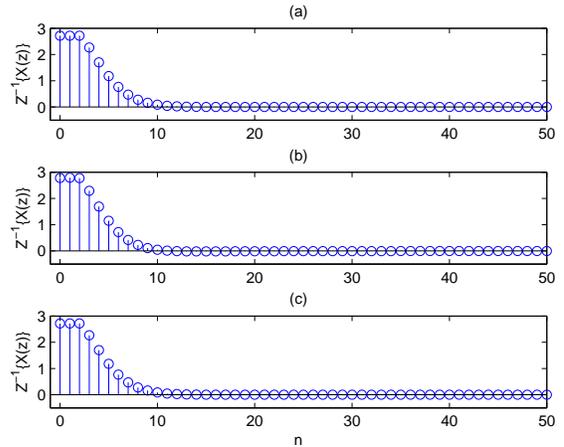}
\caption{Inverse Z-transform of $X(z)$ when (a) the first method,
(b) the second method, and (c) Eq. (\ref{int_four}) is applied
(corresponding to Example 2).}\label{fig_ex21}
\end{center}
\end{figure}

\begin{figure}[tb]
\begin{center}
\includegraphics[width=8.5cm]{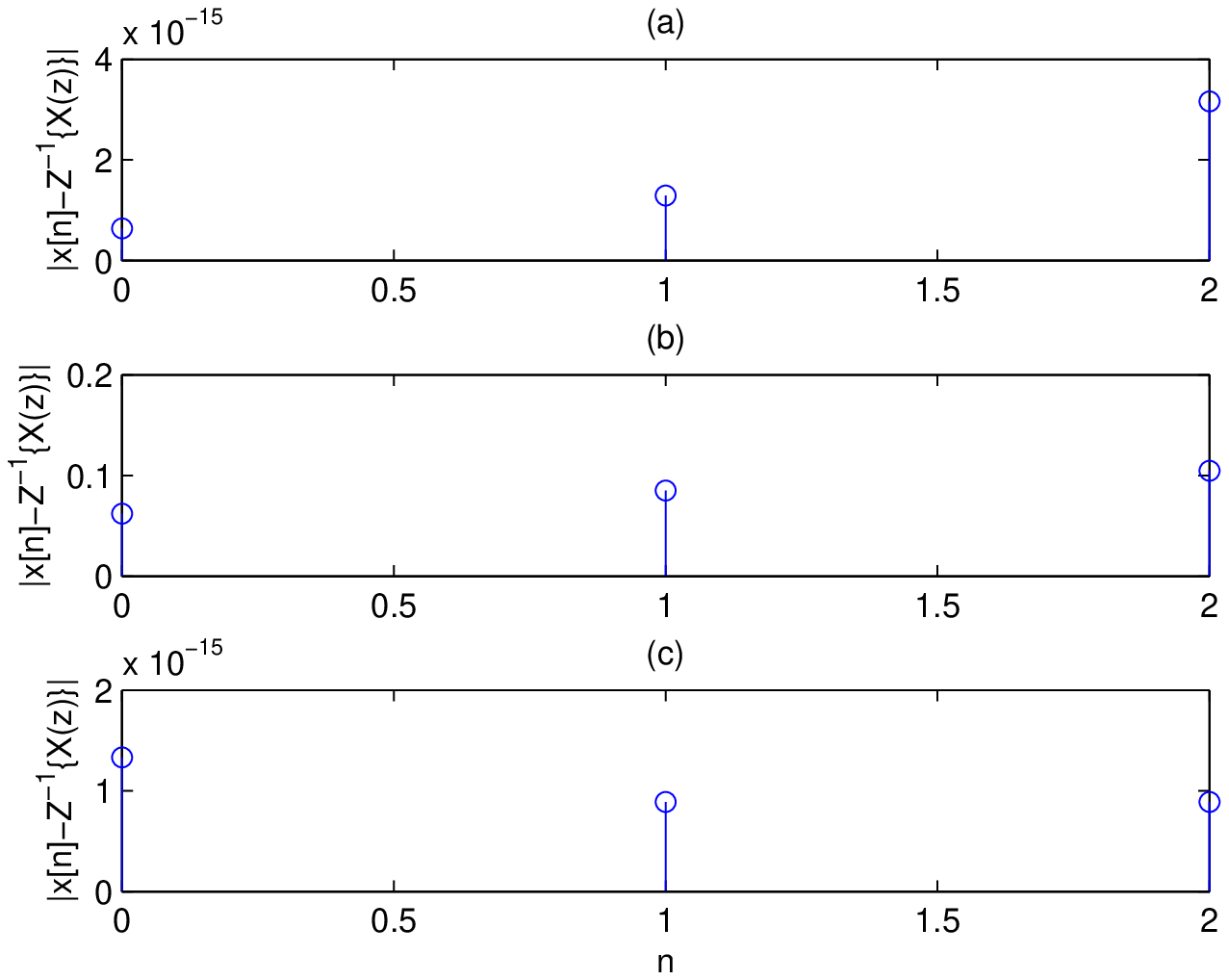}
\caption{The absolute error caused by (a) the first method, (b)
the second method, and (c) Eq. (\ref{int_four}), in the fist three
samples of $x[n]$ (corresponding to Example 2).}\label{fig_ex22}
\end{center}
\end{figure}

\begin{figure}[tb]
\begin{center}
\includegraphics[width=8.5cm]{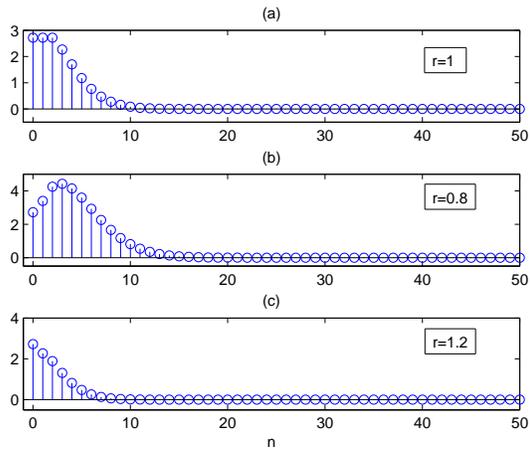}
\caption{The $x[n]$ obtained through (\ref{int_z}) when $c$ is
considered as a circle with radius $r=1,0.8,1.2$ (corresponding to
Example 2). $r=1$ leads to the most accurate
results.}\label{fig_ex23}
\end{center}
\end{figure}

\textbf{Example 2.} Consider the Z-transform given in the
following expression
\begin{equation}
X(z)=e^{e^{z^{-1}}},\quad \mathrm{ROC}:|z|>0,
\end{equation}
which can further be written as
\begin{align*}
X(z)&=e^{1+z^{-1}+z^{-2}/2!+\ldots}\\
&=e^1\times e^{z^{-1}}\times e^{z^{-2}/2!}\times \ldots\\
&=e^1\times \left(1+\frac{1}{z}+\frac{1}{2!z^2}+\ldots\right)\\
&\times \left(1+\frac{1}{2!z^2} +\frac{1}{2!
(2!z^2)^2}+\ldots\right)\times \ldots
\end{align*}
As it is observed, according to the complexity of the resulted
expression it is really difficult to extract the complete inverse
Z-transform from it (however, it can be easily verified that
$x[0]=x[1]=e$ and $x[2]=e(1/2!+1/2!)=e$). Figs.
\ref{fig_ex21}(a)-(c) show the inverse Z-transform of $X(z)$ when
the first method, second method, and Eq. (\ref{int_four}) is
applied, respectively. Figs. \ref{fig_ex22}(a)-(c) show the
absolute error caused by each of these methods. Figs.
\ref{fig_ex23}(a)-(c) show the $x[n]$ obtained through
(\ref{int_z}) when $c$ is considered as a circle with radius
$r=1,0.8,1.2$, respectively. This figure clearly shows that the
most accurate result corresponds to $r=1$.\\

\begin{thebibliography}{00}


\bibitem{pop}
A. Papoulis, Numerical inversion of the z-transform , \emph{IEEE
Transactions on Circuit Theory}, Vol. 20, No. 4, July 1973, pp.
419-420.


\bibitem{oppenheim2}
A.V. Oppenheim and R.W. Schaffer, {\it Discrete-Time Signal
Processing}, 3rd ed., Prentice Hall, NJ, 2010.

\end{thebibliography}


\end{document}